\documentclass[12pt]{elsarticle}
\usepackage{amsfonts}
\usepackage{amsmath}
\usepackage{graphicx}
\usepackage{epsfig}
\usepackage{amssymb}
\usepackage{amstext}

\usepackage{tkz-berge}

\newtheorem{theorem}{Theorem}[section]

\newenvironment{PrfFact}{{\bf Proof }}{{\hfill\tiny{$\blacksquare$\\}}}

\usepackage{tikz}
\usetikzlibrary{positioning}
\usetikzlibrary{shapes}

\usepackage{comment}

\usepackage{rotating}

\usetikzlibrary{shapes.geometric}
 				
\tikzset{
  hexagon/.style={
    shape=regular polygon, regular polygon sides=6,
    outer sep=0, minimum size=3cm, rotate=1*360/12,
    draw,
  }
}

\def\rot#1{{\rotatebox{30}{$#1$}}}
\def\rota{\rot\alpha}
\def\rotb{\rot\beta}

\newlength\flitemwidth

\begin{document}

\begin{frontmatter}

\title {An extension of a cubic $2$-connected plane graph $G$ to a  hamiltonian  plane graph contained in $G^{2}$}

\author{Jan~Florek}
\ead{jan.florek@wr.edu.pl}

\address{Faculty of Pure and Applied Mathematics,
 Wroclaw University of Science and Technology,
 Wybrze\.{z}e Wyspia\'nskiego 27,
50--370 Wroc{\l}aw, Poland}

\begin{abstract}
Let $G$ be a simple cubic $2$-connected plane graph.  For every $2$-factor $X$ of~$G$ having $n$-components  there exists a simple hamiltonian plane graph $J \subset G^{2}$ such that $|E(J)|= |E(G)| + 2n -2$ and $\Delta(J) \leqslant 5$. 

 \end{abstract}
\end{frontmatter}

\section{Introduction}

We use \cite{flobar1} as a reference for undefined terms. Let $F$ be a simple connected graph. $V(F)$  is the vertex set and $E(F)$ is the edge set of  $F$.  A \textsl{spanning subgraph} of $F$ is a subgraph whose vertex set is the entire vertex set of $F$. A spanning cycle (path) is called  \textsl{Hamilton cycle} (\textsl{Hamilton path}) and a spanning $k$-regular subgraph is called \textsl{$k$-factor}.  A graph is \textsl{hamiltonian} if it admits a Hamilton cycle. A graph is \textsl{hamiltonian-connected} if for every pair $u$, $v$ of distinct vertices of $F$, there exists a Hamilton $u-v$ path.

Given a positive integer $k$, we denote by $F^k$ the graph on $V(F)$ in which two vertices are adjacent if and only if they have distance at most $k$ in $F$. The graph $F^2$ and $F^3$ are also referred to as the \textsl{square} and \textsl{cube}, respectively, of $F$. Karaganis \cite{flobar4} and Sekanina \cite{flobar7} proved that the cube of a connected graph is hamiltonian-connected, and Fleischner \cite{flobar3} discovered that the square of a $2$-connected graph is hamiltonian (see also \v{R}\'{\i}ha \cite{flobar6}). The strengthened result (employing Fleischner's work) that the square of such a graph is hamiltonian-connected was proved by  Chartrand, Kappor,  and Nash-Williams \cite{flobar2}.

Let ${\cal G}$ be the family of all simple cubic $2$-connected plane graphs. We prove  the following  theorem.

\begin{theorem}\label{theorem1.1}
Let $G \in \cal G$. For every $2$-factor $X$ of $G$ having $n$-components  there exists a simple plane graph $J$, $G \subseteq J \subset G^{2}$, having a Hamilton cycle omitting all edges of $E(G)\backslash E(X)$, $|E(J)|= |E(G)| + 2n -2$ and $\Delta(J) \leqslant 5$. 
\end{theorem}
Notice that Petersen \cite{flobar5} proved that every simple bridgeless cubic graph  has a $1$-factor (see Distel \cite{flobar1} Corollary 2.2.2). It follows that that each graph of~$\cal G$ has a $2$-factor.

\section{Main result}

Let $G \in \cal G$. Each face  $f \in F(G)$  is bounded by a cycle $\partial(f)$ called a \textsl{facial cycle} of this face.
A cyclic sequence  $f_{1}f_{2} \ldots f_{k}f_{1}$ (a sequence $f_{1}f_{2} \ldots f_{k}$) of different faces  in $G$ is called a \textit{cyclic sequence of  faces} (a \textit{sequence of  faces} from $f_1$ to $f_k$, respectively)
if any two  successive faces $f_{i}, f_{i+1}$ are adjacent.  We say that an edge \textit{belongs} to a cyclic sequence of faces (to a sequence of faces) if it is a common edge of two successive faces  of this sequence. Then, we also say that this cyclic sequence of faces (sequence of faces) \textit{contains} this edge. Notice that  $B \subseteq E(G)$ is a bond of $G$ if and only if it is the set of all edges belonging to a some cyclic sequence of faces.

\begin{PrfFact}{\bf of Theorem \ref{theorem1.1}.}
Let $X$ be a $2$-factor of $G$ which has $n$ components. Without loss of generality we can assume that $n > 1$. We will define a face $2$-colouring ${a: F(G) \rightarrow \{\alpha, \beta\}}$. Fix a face $f$ of $F(G)$. For every $g \in F(G)$, $g \neq f$, there exists a sequence of faces from $f$ to $g$ containing only edges of $E(X)$, because $G$ is cubic and $X$ is a $2$-factor of $G$. We set $a(f)  = \alpha$ and we colour faces of this sequence  with  $\alpha$ and $\beta$ alternately. The colouring of $g$ is independent on the choice of the sequence of faces. Indeed, if $D$ is a cyclic sequence of faces and $B \subseteq E(X)$ is the set of all edges belonging to $D$, then $|B|$ is even, because $X$ is a $2$-factor of $G$. It follows that
\begin {enumerate}
\item [($1$)]
any two adjacent faces in $G$ are incident with the same edge belonging to $E(G) \backslash E(X)$ if and only if they are coloured identically by $a$.
\end {enumerate}
Further, every  facial cycle of $G$ has an \textit{orientation assigned} by $a$. Namely,  we can assume that  a facial cycle of an inner (the outer) face has the clockwise-orientation (counter clockwise-orientation, respectively) if and only if this face is coloured $\alpha$.

We can enumerate, by induction,  all components of $X$ as $C_{1}, C_{2},\ldots, C_{n}$ in such a way that there exists an edge $b_{i} \in E(G) \backslash E(X)$ connecting a subgraph $C_1 \cup \ldots \cup C_{i}$ with $C_{i+1}$, for every $1 \leqslant i < n$. Certainly any two edges of $M = \{b_{1}, \ldots, b_{n}\}$ are different and non-adjacent, because $G$ is cubic and $X$ is a $2$-factor. Since each $b_{i} \in M$ is connecting two different components of~$X$, we obtain:
\begin {enumerate}
\item [($2$)] each facial $3$-cycle of $G$ does not contain any edge of $M$,

\item [($3$)] each facial $4$-cycle  of $G$ contains at most one edge of $M$.
\end {enumerate}

Suppose that  $\cal K$ is the family of all faces in $G$ each of which is  incident with an edge of $M$. Let $f$ be any face of $\cal K$ and suppose that $\partial (f) = c_{1}c_{2} \ldots c_{n}c_{1}$ is  the facial cycle of $f$ with the  orientation assigned by $a$. Let $c_{i_{1}}c_{i_{1}+1}$, $c_{i_{2}}c_{i_{2}+1}\ldots, c_{i_{p}}c_{i_{p}+1}$ (where $p$  depends on the face $f$) be all successive edges of  $\partial (f)$ belonging to $M$.   By condition  $(2)$,  vertices $c_{i_{j}}$ and $c_{i_{j}+2}$ are not adjacent in $G$, for every $j = 1, \ldots, p$. Since edges of $M$ are not adjacent, we can draw new edges
$$e_{1}(f) = c_{i_{1}}c_{i_{1}+2}, e_{2}(f) =  c_{i_{2}}c_{i_{2}+2},\ldots, e_{p}(f) = c_{i_{p}}c_{i_{p}+2}$$
 in such a way that they are not crossing and their interiors are contained in~$f$. By adding to $G$ all edges of $\bigcup_{f \in \cal K}  \{e_{1}(f), \ldots, e_{p}(f)\}$ we obtain a plane graph~$J$ such that $G \subseteq J \subset G^{2}$ and $|E(J)|= |E(G)| + 2n -2$ . Since $G$ is simple, by conditions $(2)-(3)$, $J$ is simple too.

Let  $e =ab$ be any edge of $M$ and suppose that $e$ is incident with faces $f_{1}, f_{2} \in F(G)$. Since $e \in E(G) \backslash E(X)$, by condition $(1)$,  faces $f_1$ and $f_2$ are coloured the same by $a$. 
Assume that $a, b, c$ and $b, a, d$ are successive vertices of $\partial(f_1)$ and $\partial(f_2)$, respectively. If faces $f_1, f_2$ are coloured $\alpha$ (or $\beta$), then a $4$-cycle $D_e = adbca$ is called a \textit{diamond of type $\alpha$} (\textit{diamond of type $\beta$},  respecticely). Let  $D^{1}_{e}$ (or $D^{2}_{e}$) denote a set of two edges of $D_{e}$ belonging to $E(G)$ (to $E(J)\backslash E(G)$, respectively).  Then, $D^{1}_{e} = \{ad, bc\}$ and $D^{2}_{e} = \{db, ca\}$ (see Fig 1 and Fig 2).
Since edges of $M$ are not adjacent, from condition $(3)$ it follows that 
\begin {enumerate}
\item [($4$)]
any two different diamonds have no common edge.
\end{enumerate}
Since $G$ is cubic each vertex of $G$ belongs to at most two  diamonds. Hence, $\Delta(J) \leqslant 5$. 

\begin{figure}
\centering
\begin{tikzpicture}
  \begin{scope}[
every node/.style={anchor=west,
regular polygon,
regular polygon sides=6,
draw,
minimum width=2cm,
outer sep=0,
  rotate=-30 },
      transform shape]

    \node (A) {\rotatebox{30}{$\alpha$}};
    \node (B) at (A.corner 1) {\rotb};
    \node (C) at (B.corner 5) {\rota};
    \node (D) at (A.corner 5) {\rota};

 \filldraw[fill=lightgray]   (A.corner 4) -- (A.corner 6) -- (D.corner 1)  -- (D.corner 3);

 \draw (A.corner 6)   -- (D.corner 3);

  \filldraw[fill=lightgray]   (C.corner 3) -- (C.corner 2) -- (C.corner 4)  -- (D.corner 5);

  \draw (C.corner 3)   -- (C.corner 4);
  \draw (C.corner 3)   -- (D.corner 5);

\draw[very thick] (A.corner 3) -- (A.corner 2) -- (A.corner 1) -- (A.corner 6) -- (B.corner 5) -- (B.corner 6) -- (C.corner 1) -- (C.corner 6);

\draw[very thick] (B.corner 2) -- (B.corner 1);

\draw[very thick] (A.corner 4) -- (A.corner 5) -- (D.corner 4);

\draw[very thick]  (D.corner 5)  -- (D.corner 6) -- (C.corner 5);
\end{scope}

\node[draw=none,xshift=-1mm,yshift=-2mm] at (A.corner 4) {$d$};
\node[draw=none,xshift=-2mm,yshift=-2mm] at (A.corner 5) {$a$};

\node[draw=none,xshift=2mm,yshift=1mm] at (B.corner 4) {$b$};
\node[draw=none,xshift=0mm,yshift=2mm] at (B.corner 5) {$c$};
\node[draw=none,xshift=2mm,yshift=1mm] at (B.corner 6) {$g$};

\node[draw=none,xshift=1mm,yshift=-2mm] at (D.corner 5) {$f$};
\node[draw=none,xshift=0mm,yshift=-2mm] at (D.corner 6) {$e$};

\end{tikzpicture}

\caption{A subgraph of $G$ (without edges $bd, ac, cf, eg$) and a subgraph of a $2$-factor $X$ of $G$ (in bold). Diamonds $D(ab) = adbca$ and $D(ce) = cgefc$ are of type $\alpha$.}
\end{figure}

\begin{figure}
\centering
\begin{tikzpicture}
  \begin{scope}[%
every node/.style={anchor=west,
regular polygon,
regular polygon sides=6,
draw,
minimum width=2cm,
outer sep=0,
  rotate=-30 },
      transform shape]

    \node (A) {\rotatebox{30}{$\alpha$}};
    \node (B) at (A.corner 1) {\rotb};
    \node (C) at (B.corner 5) {\rotb};
    \node (D) at (A.corner 5) {\rota};

 \filldraw[fill=lightgray]   (A.corner 4) -- (A.corner 6) -- (D.corner 1)  -- (D.corner 3);

 \draw (A.corner 6)   -- (D.corner 3);

  \filldraw[fill=lightgray]   (B.corner 5) -- (B.corner 1) -- (B.corner 6)  -- (C.corner 4) -- cycle;

 \draw (C.corner 2)   -- (C.corner 3);

\draw[very thick] (A.corner 2) -- (A.corner 3) -- (A.corner 4) -- (A.corner 5) -- (D.corner 4) -- (D.corner 5)  ;

\draw[very thick] (B.corner 1) -- (B.corner 6) -- (C.corner 1);

\draw[very thick]  (B.corner 2) -- (B.corner 3) -- (B.corner 4) --(B.corner 5)  -- (C.corner 3) -- (C.corner 4) -- (C.corner 5) -- (C.corner 6) ;
\end{scope}

\node[draw=none,xshift=-1mm,yshift=-2mm] at (A.corner 4) {$d$};
\node[draw=none,xshift=-2mm,yshift=-2mm] at (A.corner 5) {$a$};

\node[draw=none,xshift=2mm,yshift=1mm] at (B.corner 4) {$b$};
\node[draw=none,xshift=-2mm,yshift=-2mm] at (B.corner 5) {$c$};
\node[draw=none,xshift=2mm,yshift=1mm] at (B.corner 6) {$e$};
\node[draw=none,xshift=2mm,yshift=1mm] at (B.corner 1) {$g$};

\node[draw=none,xshift=0mm,yshift=-2.5mm] at (D.corner 6) {$f$};

\end{tikzpicture}

\caption{A subgraph of $G$ (without edges $bd, ac, cg, ef$) and a subgraph of a $2$-factor $X$ of $G$ (in bold). Diamons $D(ab) = adbca$ is of type $\alpha$ and $D(ce)=cgefc$ is of type $\beta$.}
\end{figure}
Define, by induction,  a  subgraph of $J$: $H_1 = C_1$ and
$$H_{i+1} = (C_{1} \cup \ldots \cup C_{i+1} - (D^{1}_{b_1}\cup \ldots \cup D^{1}_{b_i})) +  (D^{2}_{b_1}\cup \ldots \cup D^{2}_{b_i}), \hbox{ for } 1 \leqslant i < n.$$
Notice that $H_{i}$ is  omitting all edges of $E(G)\backslash E(X)$, because $C_{j}$ is contained in $X$, for $1 \leqslant j \leqslant n$, and $D^{2}_{b_j}$ is contained in $E(J) \backslash E(G)$, for $1 \leqslant j < n$. Further, $H_{i}$ is disjoint with $C_{i+1}$ because $C_{i+1}$ is disjoint with $C_j$, for $j < i$.

Assume that $H_i$ is  a cycle of $J$. We prove that $H_{i+1}$ is also a cycle of $J$.  Since $b_i$ is an edge connecting $C_{1} \cup \ldots \cup C_{i}$ with $C_{i+1}$, by condition $(4)$, any edge of $D^{1}_{b_i}$ is not an edge of  $D^{1}_{b_j}$, for $j < i$. Hence, one edge of  $D^{1}_{b_i}$ is an edge of $H_{i+1}$ and another one is an edge of $C_{i+1}$. Further,  by condition~$(4)$, any edge of $D^{2}_{b_i}$ is not an edge of  $D^{2}_{b_j}$, for $j < i$. Hence, any edge of  $D^{2}_{b_i}$ is not an edge of $H_{i} \cup C_{i+1}$.  Therefore, $H_{i+1}$ is a cycle  of $J$. Thus, $H_n$  is a Hamilton cycle of $J$, because $V(H_n) = V(C_{1}) \cup\ldots \cup V(C_{n})= V(G) = V(J)$. 
\end{PrfFact}


\end{document}